# UTILITY MAXIMIZATION IN INCOMPLETE MARKETS[1]

By Ying Hu, Peter Imkeller and Matthias Müller

*Université de Rennes 1, Humboldt-Universität zu Berlin and Humboldt-Universität zu Berlin*

We consider the problem of utility maximization for small traders on incomplete financial markets. As opposed to most of the papers dealing with this subject, the investors' trading strategies we allow underly constraints described by closed, but not necessarily convex, sets. The final wealths obtained by trading under these constraints are identified as stochastic processes which usually are supermartingales, and even martingales for particular strategies. These strategies are seen to be optimal, and the corresponding value functions determined simply by the initial values of the supermartingales. We separately treat the cases of exponential, power and logarithmic utility.

**Introduction.** In this paper we consider a small trader on an incomplete financial market who can trade in a finite time interval $[0,T]$ by investing in risky stocks and a riskless bond. He aims at maximizing the utility he draws from his final wealth measured by some utility function. The trading strategies he may choose to attain his wealth underly some restriction formalized by a constraint. For example, he may be forced not to have a negative number of shares or that his investment in risky stocks is not allowed to exceed a certain threshold. We will be interested not only in describing the trader's optimal utility, but also the strategies which he may follow to reach this goal. As opposed to most of the papers dealing so far with the maximization of expected utility under constraints, we essentially relax the hypotheses to be fulfilled by them. They are formulated as usual by the requirement that the strategies take their values in some set, which is supposed to simply

Received December 2003; revised October 2004.

[1] Supported in part by the DFG research center "Mathematics for key technologies" (FZT 86) in Berlin.

be closed instead of convex. We consider three types of utility functions. In Section 2 we carry out the calculation of the value function and an optimal strategy for exponential utility. In this case, the investor is allowed to have an additional liability, and maximizes the utility of its sum with terminal wealth. In Section 3 we consider power utility, and in Section 4 the simplest one: logarithmic utility.

The method that we apply in order to obtain value function and optimal strategy is simple. We propose to construct a stochastic process $R^\rho$ depending on the investor's trading strategy $\rho$, and such that its terminal value equals the utility of the trader's terminal wealth. As mentioned above, to model the constraint, trading strategies are supposed to take their values in a closed set. In our market, the absence of completeness is not explicitly described by a set of martingale measures equivalent to the historical probability. Instead, we choose $R^\rho$ such that for every trading strategy $\rho$, $R^\rho$ is a supermartingale. Moreover, there exists at least one particular trading strategy $\rho^*$ such that $R^{\rho^*}$ is a martingale. Hereby, the initial value is supposed to not depend on the strategy. Evidently, the strategy $\rho^*$ related to the martingale has to be the optimal one. Then the value function of the optimization problem is just given by the initial value of $R^{\rho^*}$.

Since we work on a Wiener filtration, the powerful tool of backward stochastic differential equations (BSDE) is available. It allows the construction of the stochastic control process $\rho^*$, and thus the description of the value function in terms of the solution of a BSDE.

In a related paper, El Karoui and Rouge [7] compute the value function and the optimal strategy for exponential utility by means of BSDE, assuming more restrictively that the strategies be confined to a convex cone. Sekine [15] relies on a duality result obtained by Cvitanic and Karatzas [2], also describing constraints through convex cones. He studies the maximization problem for the exponential and power utility functions, and uses an attainability condition which solves the primal and dual problems, finally writing this condition as a BSDE. In contrast to these papers, we do not use duality, and directly characterize the solution of the primal problem. This allows us to pass from convex to closed constraints.

Utility maximization is one of the most frequent problems in financial mathematics and has been considered by numerous authors. Here are some of the milestones viewed from our perspective of maximization under constraints using the tools of BSDEs. For a complete market, utility maximization has been considered in [9]. Cvitanic and Karatzas [2] prove existence and uniqueness of the solution for the utility maximization problem in a Brownian filtration constraining strategies to convex sets. There are numerous papers considering general semimartingales as stock price processes. Delbaen et al. [4] give a duality result between the optimal strategy for the maximization of the exponential utility and the martingale measure minimizing the relative entropy with respect to the real world measure $P$. This duality can be used to characterize the utility indifference price for an option. Also relying upon duality theory, Kramkov and Schachermayer [12] and Cvitanic, Schachermayer and Wang [3] give a fairly complete solution of the utility optimization problem on incomplete markets for a class of general utility functions not containing the exponential one. See also the review paper by Schachermayer [16] for a more complete account and further references.

The powerful tool of BSDE has been introduced to stochastic control theory by Bismut [1]. Its mathematical treatment in terms of stochastic analysis was initiated by Pardoux and Peng [14], and its particular significance for the field of utility maximization in financial stochastics clarified in [6].

**1. Preliminaries and the market model.** A probability space $(\Omega, \mathcal{F}, P)$ carrying an $m$-dimensional Brownian motion $(W_t)_{t \in [0,T]}$ is given. The filtration $\mathbb{F}$ is the completion of the filtration generated by $W$.

Let us briefly explain some special notation that will be used in the paper. $|\cdot|$ stands for the Euclidean norm in $\mathbb{R}^m$. For $q \geq 1$, $L^q$ denotes the set of $\mathcal{F}_T$-measurable random variables $F$ such that $E[|F|^q] < \infty$, for $k \in \mathbb{N}$, $\mathcal{H}^k(\mathbb{R}^d)$ the set of all $\mathbb{R}^d$-valued stochastic processes $\vartheta$ which are predictable with respect to $\mathbb{F}$ and satisfy $E[\int_0^T |\vartheta_t|^k \, dt] < \infty$. $\mathcal{H}^\infty(\mathbb{R}^d)$ is the set of all $\mathbb{F}$-predictable $\mathbb{R}^d$-valued processes that are $\lambda \otimes P$-a.e. bounded on $[0,T] \times \Omega$. Note here that we write $\lambda$ for the Lebesgue measure on $[0,T]$ or $\mathbb{R}$.

Let $M$ denote a continuous semimartingale. The stochastic exponential $\mathcal{E}(M)$ is given by

$$\mathcal{E}(M)_t = \exp(M_t - \tfrac{1}{2}\langle M \rangle_t), \qquad t \in [0,T],$$

where the quadratic variation is denoted by $\langle M \rangle$. Let $C$ denote a closed subset of $\mathbb{R}^m$ and $a \in \mathbb{R}^m$. The distance between $a$ and $C$ is defined as

$$\text{dist}_C(a) = \min_{b \in C} |a - b|.$$

The set $\Pi_C(a)$ consists of those elements of $C$ at which the minimum is obtained:

$$\Pi_C(a) = \{b \in C : |a - b| = \text{dist}_C(a)\}.$$

This set is not empty and evidently may contain more than one point.

The financial market consists of one bond with interest rate zero and $d \leq m$ stocks. In case $d < m$, we face an incomplete market. The price process of stock $i$ evolves according to the equation

(1) $$\frac{dS_t^i}{S_t^i} = b_t^i \, dt + \sigma_t^i \, dW_t, \qquad i = 1, \ldots, d,$$

where $b^i$ (resp. $\sigma^i$) is an $\mathbb{R}$-valued (resp. $\mathbb{R}^{1\times m}$-valued) predictable uniformly bounded stochastic process. The lines of the $d \times m$-matrix $\sigma$ are given by the vector $\sigma_t^i$, $i=1,\ldots,d$. The volatility matrix $\sigma = (\sigma^i)_{i=1,\ldots,d}$ has full rank and we assume that $\sigma\sigma^{tr}$ is uniformly elliptic, that is, $KI_d \geq \sigma\sigma^{tr} \geq \varepsilon I_d$, $P$-a.s. for constants $K > \varepsilon > 0$. The predictable $\mathbb{R}^m$-valued process

$$\theta_t = \sigma_t^{tr}(\sigma_t\sigma_t^{tr})^{-1}b_t, \qquad t \in [0,T],$$

is then also uniformly bounded.

A $d$-dimensional $\mathbb{F}$-predictable process $\pi = (\pi_t)_{0 \leq t \leq T}$ is called trading strategy if $\int \pi \frac{dS}{S}$ is well defined, for example, $\int_0^T \|\pi_t\sigma_t\|^2\,dt < \infty$ $P$-a.s. For $1 \leq i \leq d$, the process $\pi_t^i$ describes the amount of money invested in stock $i$ at time $t$. The number of shares is $\frac{\pi_t^i}{S_t^i}$. The wealth process $X^\pi$ of a trading strategy $\pi$ with initial capital $x$ satisfies the equation

$$X_t^\pi = x + \sum_{i=1}^d \int_0^t \frac{\pi_{i,u}}{S_{i,u}}\,dS_{i,u} = x + \int_0^t \pi_u\sigma_u(dW_u + \theta_u\,du), \qquad t \in [0,T].$$

In this notation $\pi$ has to be taken as a vector in $\mathbb{R}^{1\times d}$. Trading strategies are self-financing. The investor uses his initial capital and during the trading interval $[0,T]$, there is no extra money flow out of or into his portfolio. Gains or losses are only obtained by trading with the stock.

The optimal trading strategy we will find in this paper happens to be in the class of martingales of bounded mean oscillation, briefly called BMO-martingales. Here we recall a few well-known facts from this theory following the exposition in [10]. The statements in [10] are made for infinite time horizon. In the text they will be applied to the simpler framework of finite time horizon, replacing $\infty$ with $T$. Let $\mathbb{G}$ be a complete, right-continuous filtration, $P$ a probability measure and $M$ a continuous local $(P,\mathbb{G})$-martingale satisfying $M_0 = 0$. Let $1 \leq p < \infty$. Then $M$ is in the normed linear space $\mathrm{BMO}_p$ if

$$\|M\|_{\mathrm{BMO}_p} := \sup_{\tau\ \mathbb{G}\text{-stopping time}} E[|M_T - M_\tau|^p|\mathcal{G}_\tau]^{1/p} < \infty.$$

By Corollary 2.1 in [10], $M$ is a $\mathrm{BMO}_p$-martingale if and only if it is a $\mathrm{BMO}_q$-martingale for every $q \geq 1$. Therefore, it is simply called a BMO-martingale. In particular, $M$ is a BMO-martingale if and only if

$$\|M\|_{\mathrm{BMO}_2} = \sup_{\tau\ \mathbb{G}\text{-stopping time}} E[\langle M\rangle_T - \langle M\rangle_\tau|\mathcal{G}_\tau]^{1/2} < \infty.$$

This means local martingales of the form $M_t = \int_0^t \xi_s\,dW_s$ are BMO-martingales if and only if

$$(2) \qquad \|M\|_{\mathrm{BMO}_2} = \sup_{\tau\ \mathbb{G}\text{-stopping time}} E\left[\int_\tau^T \|\xi_s\|^2\,ds\Big|\mathcal{G}_\tau\right]^{1/2} < \infty.$$

Due to the finite time horizon, this condition is satisfied for bounded integrands. According to Theorem 2.3 in [10], the stochastic exponential $\mathcal{E}(M)$ of a BMO-martingale $M$ is a uniformly integrable martingale. If $Q$ is a probability measure defined by $dQ = \mathcal{E}(M)_T \, dP$ for a $P$-BMO martingale $M$, then the Girsanov transform of a $P$-BMO martingale is a BMO-martingale under $Q$ (Theorem 3.6 in [10]).

Suppose our investor has a liability $F$ at time $T$. This random variable $F$ is assumed to be $\mathcal{F}_T$-measurable and bounded, but not necessarily positive. He tries to find a trading strategy that is optimal in presence of this liability $F$, in a sense to be made precise in the beginning of the following section.

In order to compute the optimal trading strategy, we use quadratic Backward Stochastic Differential Equations (BSDE) and apply a result of Kobylanski [11] to get existence of a solution for our BSDE. This result is proved for bounded terminal random variables. Therefore, we have to assume that $F$ is bounded.

**2. Exponential utility.** In this section we specify the sense of optimality for trading strategies by stipulating that the investor wants to maximize his expected utility with respect to the exponential utility from his total wealth $X_T^p - F$. Let us recall that, for $\alpha > 0$, the exponential utility function is defined as

$$U(x) = -\exp(-\alpha x), \qquad x \in \mathbb{R}.$$

The definition of admissible trading strategies guarantees that there is no arbitrage. In addition, we allow constraints on the trading strategies. Formally, they are supposed to take their values in a closed set, that is, $\pi_t(\omega) \in \tilde{C}$, with $\tilde{C} \subseteq \mathbb{R}^{1 \times d}$. We emphasize that $\tilde{C}$ is not assumed to be convex.

DEFINITION 1 (*Admissible strategies with constraints*). Let $\tilde{C}$ be a closed set in $\mathbb{R}^{1 \times d}$. The set of admissible trading strategies $\tilde{\mathcal{A}}$ consists of all $d$-dimensional predictable processes $\pi = (\pi_t)_{0 \leq t \leq T}$ which satisfy $E[\int_0^T |\pi_t \sigma_t|^2 \, dt] < \infty$ and $\pi_t \in \tilde{C}$ $\lambda \otimes P$-a.s., as well as

$$\{\exp(-\alpha X_\tau^\pi) : \tau \text{ stopping time with values in } [0, T]\}$$

is a uniformly integrable family.

REMARK 2. The condition of square integrability in Definition 1 guarantees that there is no arbitrage. In fact, the square integrability condition on $\pi$ and the boundedness of $\theta$ yields that $E[\sup_{0 \leq t \leq T}(X_t^\pi)^2] < \infty$. According to Theorem 2.1 in [14], $(X_t, \pi_t \sigma_t)$ is the unique solution of the BSDE

$$X_t = X_T - \int_t^T (\pi_s \sigma_s) \, dW_s - \int_t^T (\pi_s \sigma_s) \theta_s \, ds,$$

with $E[\int_0^T (X_s^\pi)^2 \, ds] < \infty$, $E[\int_0^T (\pi_s \sigma_s)^2 \, ds] < \infty$. So the initial capital $X_0^\pi$ needed to attain $X_T^\pi$ is uniquely determined. In particular, Theorem 2.2 in [6] yields if $X_0^\pi = 0$ and $X_T^\pi \geq 0$ $P$-a.s., then $X_T^\pi = 0$ $P$-a.s.

REMARK 3. In accordance with the classical literature (see [5]) the uniform integrability condition in Definition 1 coincides with the notion of class D.

REMARK 4. If $X^\pi$ is square integrable and $\pi_t \in \tilde{C}$ $\lambda \otimes P$-a.s., as well as $X^\pi$ is bounded from below on $[0,T]$, it is obvious that $\pi \in \tilde{\mathcal{A}}$.

For $t \in [0,T], \omega \in \Omega$ define the set $C_t(\omega) \subseteq \mathbb{R}^m$ by

$$C_t(\omega) = \tilde{C}\sigma_t(\omega). \tag{3}$$

The entries of the matrix-valued process $\sigma$ are uniformly bounded. Therefore, we get

$$\min\{|a| : a \in C_t(\omega)\} \leq k_1 \quad \text{for } \lambda \otimes P\text{-a.e.}(t,\omega) \tag{4}$$

with a constant $k_1 \geq 0$. Furthermore, for every $(\omega, t)$, the set $C_t(\omega)$ is closed. This is crucial for our analysis.

REMARK 5. Writing

$$p_t = \pi_t \sigma_t, \qquad t \in [0,T],$$

the set of admissible trading strategies $\tilde{\mathcal{A}}$ is equivalent to a set $\mathcal{A}$ of $\mathbb{R}^{1 \times m}$-valued predictable stochastic processes $p$ with $p \in \mathcal{A}$ iff $E[\int_0^T |p(t)|^2 \, dt] < \infty$ and $p_t(\omega) \in C_t(\omega)$ $P$-a.s., as well as

$$\{\exp(-\alpha X_\tau^p) : \tau \text{ stopping time with values in } [0,T]\}$$

is a uniformly integrable family.

Such a process $p \in \mathcal{A}$ will also be named strategy, and $X^{(p)}$ denotes its wealth process.

So the investor wants to solve the maximization problem

$$V(x) := \sup_{\pi \in \tilde{\mathcal{A}}} E\left[-\exp\left(-\alpha\left(x + \int_0^T \pi_t \frac{dS_t}{S_t} - F\right)\right)\right],$$

where $x$ is the initial wealth. $V$ is called value function. Losses, that is, realizations with $X^\pi - F < 0$, are punished very strongly. Large gains or realizations with $X^\pi - F > 0$ are weakly valued.

REMARK 6. We shall show below that the sup is taken by a particular strategy $p^*$ which is admissible in the sense of our definition. Note that this process might not lead to a wealth process which is bounded from below, and therefore not admissible in this sense. For further details, see [13] and [17].

The maximization problem is evidently equivalent to

$$(5) \qquad V(x) = \sup_{p \in \mathcal{A}} E\bigg[-\exp\bigg(-\alpha\bigg(x + \int_0^T p_t(dW_t + \theta_t \, dt) - F\bigg)\bigg)\bigg].$$

In order to find the value function and an optimal strategy, we construct a family of stochastic processes $R^{(p)}$ with the following properties:

- $R_T^{(p)} = -\exp(-\alpha(X_T^p - F))$ for all $p \in \mathcal{A}$,
- $R_0^{(p)} = R_0$ is constant for all $p \in \mathcal{A}$,
- $R^{(p)}$ is a supermartingale for all $p \in \mathcal{A}$ and there exists a $p^* \in \mathcal{A}$ such that $R^{(p^*)}$ is a martingale.

The process $R^{(p)}$ and its initial value $R_0$ depend, of course, on the initial capital $x$. Given processes possessing these properties, we can compare the expected utilities of the strategies $p \in \mathcal{A}$ and $p^* \in \mathcal{A}$ by

$$(6) \quad E[-\exp(-\alpha(X_T^p - F))] \le R_0(x) = E[-\exp(-\alpha(X_T^{p^*} - F))] = V(x),$$

whence $p^*$ is the desired optimal strategy. To construct this family, we set

$$R_t^{(p)} := -\exp(-\alpha(X_t^{(p)} - Y_t)), \qquad t \in [0, T], p \in \mathcal{A},$$

where $(Y, Z)$ is a solution of the BSDE

$$Y_t = F - \int_t^T Z_s \, dW_s - \int_t^T f(s, Z_s) \, ds, \qquad t \in [0, T].$$

In these terms we are bound to choose a function $f$ for which $R^{(p)}$ is a supermartingale for all $p \in \mathcal{A}$ and there exists a $p^* \in \mathcal{A}$ such that $R^{(p^*)}$ is a martingale. This function $f$ also depends on the constraint set $(C_t)$, where $(p_t)$ takes its values [see (3)]. We get

$$V(x) = R_0^{(p,x)} = -\exp(-\alpha(x - Y_0)) \qquad \text{for all } p \in \mathcal{A}.$$

In order to calculate $f$, we write $R$ as the product of a (local) martingale $M^{(p)}$ and a (not strictly) decreasing process $\tilde{A}^{(p)}$ that is constant for some $p^* \in \mathcal{A}$. For $t \in [0, T]$, define

$$M_t^{(p)} = \exp(-\alpha(x - Y_0)) \exp\bigg(-\int_0^t \alpha(p_s - Z_s) \, dW_s - \tfrac{1}{2} \int_0^t \alpha^2 (p_s - Z_s)^2 \, ds\bigg).$$

Comparing $R^{(p)}$ and $M^{(p)} \tilde{A}^{(p)}$ yields

$$\tilde{A}_t^{(p)} = -\exp\left(\int_0^t v(s, p_s, Z_s)\, ds\right), \qquad t \in [0, T],$$

with

$$v(t, p, z) = -\alpha p \theta_t + \alpha f(t, z) + \tfrac{1}{2}\alpha^2 |p - z|^2.$$

In order to obtain a decreasing process $\tilde{A}^{(p)}$, evidently $f$ has to satisfy

$$v(t, p_t, Z_t) \geq 0 \qquad \text{for all } p \in \mathcal{A}$$

and

$$v(t, p_t^*, Z_t) = 0$$

for some particular $p^* \in \mathcal{A}$. For $t \in [0, T]$, we have

$$\begin{aligned}
\frac{1}{\alpha} v(t, p_t, Z_t) &= \frac{\alpha}{2}|p_t|^2 - \alpha p_t\left(Z_t + \frac{1}{\alpha}\theta_t\right) + \frac{\alpha}{2}|Z_t|^2 + f(t, Z_t) \\
&= \frac{\alpha}{2}\left|p_t - \left(Z_t + \frac{1}{\alpha}\theta_t\right)\right|^2 - \frac{\alpha}{2}\left|Z_t + \frac{1}{\alpha}\theta_t\right|^2 + \frac{\alpha}{2}Z_t^2 + f(t, Z_t) \\
&= \frac{\alpha}{2}\left|p_t - \left(Z_t + \frac{1}{\alpha}\theta_t\right)\right|^2 - Z_t \theta_t - \frac{1}{2\alpha}|\theta_t|^2 + f(t, Z_t).
\end{aligned}$$

Now set

$$f(t, z) = -\frac{\alpha}{2}\operatorname{dist}^2\left(z + \frac{1}{\alpha}\theta_t, C_t(\omega)\right) + z\theta_t + \frac{1}{2\alpha}|\theta_t|^2.$$

For this choice, we get $v(t, p, z) \geq 0$ and for

$$p_t^* \in \Pi_{C_t(\omega)}\left(Z_t + \frac{1}{\alpha}\theta_t\right), \qquad t \in [0, T],$$

we obtain $v(\cdot, p^*, Z) = 0$.

Here we see why the set $\tilde{C}$ and, hence, $C_t$ on which trading strategies are restricted is assumed to be closed. In order to find the value function, we have to minimize the distance between a point and a set. Furthermore, there must exist some element in $C_t$ realizing the minimal distance. Both requirements are satisfied for closed sets. In a convex set the minimizer is unique. This would lead to a unique utility maximizing trading strategy. However, we prove existence of a possibly nonunique trading strategy solving the maximization problem for closed but not necessarily convex constraints.

THEOREM 7. *The value function of the optimization problem* (5) *is given by*

$$V(x) = -\exp(-\alpha(x - Y_0)),$$

where $Y_0$ is defined by the unique solution $(Y, Z) \in \mathcal{H}^\infty(\mathbb{R}) \times \mathcal{H}^2(\mathbb{R}^m)$ of the BSDE

(7) $$Y_t = F - \int_t^T Z_s \, dW_s - \int_t^T f(s, Z_s) \, ds, \qquad t \in [0, T],$$

with

$$f(\cdot, z) = -\frac{\alpha}{2} \operatorname{dist}^2\left(z + \frac{1}{\alpha}\theta, C\right) + z\theta + \frac{1}{2\alpha}|\theta|^2.$$

There exists an optimal trading strategy $p^* \in \mathcal{A}$, with

(8) $$p_t^* \in \Pi_{C_t(\omega)}\left(Z_t + \frac{1}{\alpha}\theta_t\right), \qquad t \in [0, T], P\text{-a.s.}$$

PROOF. In order to get the existence of solutions of the BSDE (7), we apply Theorem 2.3 of [11]. According to Lemma 11 below, for fixed $z \in \mathbb{R}^m$, $(f(t, z))_{t \in [0, T]}$ defines a predictable process. A sufficient condition for the existence of a solution is condition (H1) in [11]: there are constants $c_0, c_1$ such that

(9) $$|f(t, z)| \leq c_0 + c_1 |z|^2 \qquad \text{for all } z \in \mathbb{R}^n \ P\text{-a.s.}$$

By means of (4), we get, for $z \in \mathbb{R}^m, t \in [0, T]$,

$$\operatorname{dist}^2\left(z + \frac{1}{\alpha}\theta_t, C_t\right) \leq 2|z|^2 + 2\left(\frac{1}{\alpha}|\theta_t| + k_1\right)^2.$$

So (9) follows from the boundedness of $\theta$. Theorem 2.3 in [11] states that the BSDE (7) possesses at least one solution $(Y, Z) \in \mathcal{H}^\infty(\mathbb{R}) \times \mathcal{H}^2(\mathbb{R}^m)$.

To prove uniqueness, suppose that solutions $(Y^1, Z^1) \in \mathcal{H}^\infty(\mathbb{R}) \times \mathcal{H}^2(\mathbb{R}^m)$, $(Y^2, Z^2) \in \mathcal{H}^\infty(\mathbb{R}) \times \mathcal{H}^2(\mathbb{R}^m)$ of the BSDE are given. Then we have

$$Y^1 - Y^2 = -\int_\cdot^T (Z^1 - Z^2) \, dW - \int_\cdot^T (f(s, Z_s^1) - f(s, Z_s^2)) \, ds.$$

Now note that, for $s \in [0, T], z^1, z^2 \in \mathbb{R}^m$, we may write

$$f(s, z^1) - f(s, z^2)$$
$$= -\frac{\alpha}{2}\left[\operatorname{dist}^2\left(z^1 + \frac{1}{\alpha}\theta_s, C_s\right) - \operatorname{dist}^2\left(z^2 + \frac{1}{\alpha}\theta_s, C_s\right)\right] + (z^1 - z^2)\theta_s.$$

Using the Lipschitz property of the distance function from a closed set, we obtain the estimate

$$|f(s, z^1) - f(s, z^2)| \leq c_1 |z^1 - z^2| + c_2(|z^1| + |z^2|)(|z^1 - z^2|)$$
$$\leq c_3(1 + |z^1| + |z^2|)|z^1 - z^2|.$$

Let us set

$$\beta(t) = \begin{cases} \dfrac{f(t, Z_t^1) - f(t, Z_t^2)}{Z_t^1 - Z_t^2}, & \text{if } Z_t^1 - Z_t^2 \neq 0, \\ 0, & \text{if } Z_t^1 - Z_t^2 = 0. \end{cases}$$

Then we obtain from the preceding estimate

$$|\beta(t)| \leq c(1 + |Z_t^1| + |Z_t^2|), \qquad t \in [0, T].$$

Moreover, from the boundedness of $Y^1$ and $Y^2$, the $P$-BMO property of $\int_0^{\cdot} Z^i(s)\,dW_s$, $i = 1, 2$, follows, see Lemma 12. This in turn entails that $\int_0^{\cdot} \beta(s)\,dW_s$ is a $P$-BMO martingale. But this allows us to give an alternative description of the difference of solutions in

$$Y^1 - Y^2 = -\int_{\cdot}^{T} (Z_s^1 - Z_s^2)\,dW_s - \int_{\cdot}^{T} \beta(s)(Z_s^1 - Z_s^2)\,ds$$

$$= -\int_{\cdot}^{T} (Z_s^1 - Z_s^2)[dW_s + \beta(s)\,ds].$$

This process is a martingale under the equivalent probability measure $Q$, which has density

$$\mathcal{E}\left(-\int_0^T \beta(t)\,dW_t\right)$$

with respect to $P$. Since $Y_T^1 = F = Y_T^2$, we therefore conclude $Y^1 = Y^2$ and $Z^1 = Z^2$, and uniqueness is established.

To find the value function of our optimization problem, we proceed with the unique solution $(Y, Z) \in \mathcal{H}^{\infty}(\mathbb{R}) \times \mathcal{H}^2(\mathbb{R}^m)$ of (7). Let $p^*$ denote the predictable process constructed in Lemma 11 for $a = Z + \frac{1}{\alpha}\theta$. Then $\tilde{A}_t^{(p^*)}(\omega) = -1$ for $\lambda \otimes P$ almost all $(t, \omega)$. By Lemma 12 below, $\int_0^{\cdot}(p_s^* - Z_s)\,dW_s$ is a $P$-BMO martingale, whence $R^{(p^*)}$ is uniformly integrable (Theorem 2.3 in [10]). Since, moreover, $Y$ is a bounded process, we obtain the uniform integrability of the family $\{\exp(-\alpha X_{\tau}^{(p^*)}) : \tau$ stopping time in $[0, T]\}$. Therefore, $p^* \in \mathcal{A}$. Hence, $R^{(p^*, x)}$ is a martingale and

$$R_0^{(p^*)} = E\left[-\exp\left(-\alpha\left(x + \int_0^T p_s^*(dW_s + \theta_s\,ds) - F\right)\right)\right]$$

$$= -\exp(-\alpha(x - Y_0)).$$

It remains to show that $R^{(p)}$ is a supermartingale for all $p \in \mathcal{A}$. Since $p \in \mathcal{A}$, the process $M = M_0 \mathcal{E}(-\alpha \int (p_s - Z_s)\,dW_s)$ is a local martingale. Hence, there exists a sequence of stopping times $(\tau_n)_{n \in \mathbb{N}}$ satisfying $\lim_{n \to \infty} \tau_n = T$ $P$-a.s. such that $(M_{t \wedge \tau_n})_t$ is a positive martingale for each $n \in \mathbb{N}$. The process $\tilde{A}^{(p)}$ is decreasing. Thus, $R^{(p)}_{t\wedge\tau_n} = M_{t\wedge\tau_n}\tilde{A}^{(p)}_{t\wedge\tau_n}$ is a supermartingale, that is, for $s \leq t$,

$$E[R^{(p)}_{t\wedge\tau_n}|\mathcal{F}_s] \leq R^{(p)}_{s\wedge\tau_n}.$$

For any set $A \in \mathcal{F}_s$, we have

$$E[R^{(p)}_{t\wedge\tau_n}\mathbb{1}_A] \leq E[R^{(p)}_{s\wedge\tau_n}\mathbb{1}_A].$$

Since $\{R^{(p)}_{t\wedge\tau_n}\}_n$ and $\{R^{(p)}_{s\wedge\tau_n}\}_n$ are uniformly integrable by the definition of admissibility and the boundedness of $Y$, we may let $n$ tend to $\infty$ to obtain

$$E[R^{(p)}_t\mathbb{1}_A] \leq E[R^{(p)}_s\mathbb{1}_A].$$

This implies the claimed supermartingale property of $R^{(p)}$. $\square$

REMARK 8. If the process $\int_0^\cdot p_s\,dW_s$ is a BMO martingale and $E[\exp\times(-\alpha(X^{(p)}_T - F))] < \infty$, a variant of an argument of the above proof can be used to see that $p \in \mathcal{A}$. In fact, we see that $M^{(p)}$ is a uniformly integrable martingale, while $A^{(p)}$ is decreasing. Hence, $R^{(p)}$ is a supermartingale. This just states that, for stopping times $\tau$,

$$-\exp(-\alpha(X^{(p)}_\tau - Y_\tau)) \geq E[-\exp(-\alpha(X^{(p)}_T - F))|\mathcal{F}_\tau].$$

Consequently,

$$\exp(-\alpha X^{(p)}_\tau) \leq \exp(-\alpha Y_\tau)E[\exp(-\alpha(X^{(p)}_T - F))|\mathcal{F}_\tau].$$

This clearly implies uniform integrability of $\{\exp(-\alpha X^{(p)}_\tau) : \tau$ stopping time in $[0, T]\}$.

We can show that the strategy $p^*$ is optimal in a wider sense. In fact, an investor who has chosen at time 0 the strategy $p^*$ will stick to this decision if he starts solving the optimization problem at some later time between 0 and $T$. For this purpose, let us formulate the optimization problem more generally for a stopping time $\tau \leq T$ and an $\mathcal{F}_\tau$-measurable random variable which describes the capital at time $\tau$, that is, $X_\tau = X^p_\tau$ for some $p \in \mathcal{A}$. So we consider the maximization problem

$$V(\tau, X_\tau) = \operatorname*{ess\,sup}_{p \in \mathcal{A}} E\left[-\exp\left(-\alpha\left(X_\tau + \int_\tau^T p_s(dW_s + \theta_s\,ds) - F\right)\right)\Big|\mathcal{F}_\tau\right].$$
(10)

PROPOSITION 9 (Dynamic principle). *The value function $x \mapsto -\exp(-\alpha(x - y))$ satisfies the dynamic programming principle, that is,*

$$V(\tau, X_\tau) = -\exp(-\alpha(X_\tau - Y_\tau))$$

for all stopping times $\tau \leq T$, where $Y_\tau$ belongs to a solution of the BSDE (7). An optimal strategy that attains the essential supremum in (10) is given by $p^*$, the optimal strategy constructed in Theorem 7.

PROOF. For $t \in [0,T]$, set

$$R_t = -\exp(-\alpha(X_t - Y_t))\mathcal{E}\left(-\int_t^T \alpha(p_s - Z_s)\,dW_s\right)\exp\left(\int_t^T v(s, p_s, Z_s)\,ds\right)$$

and apply the optional stopping theorem to the stochastic exponential. The claim follows as in Theorem 7. □

REMARK 10. If the constraint $C$ on the strategies is a convex cone, the value function $V$ and the optimal strategy $p^*$ both constructed in Theorem 7 are equivalent to those determined in [15] and [7].

Sekine considers the utility function $x \mapsto -\frac{1}{\alpha}\exp(-\alpha x)$. He obtains the value function

$$V(x) = -\frac{1}{\alpha}\exp(-\alpha x + \bar{Y}_0),$$

starting with the BSDE

$$\bar{Y}_t = \alpha F - \int_t^T \bar{z}_s\,dW_s - \int_t^T \bar{f}(s, \theta_s, \bar{z}_s)\,ds, \qquad t \in [0,T],$$

where

$$\bar{f}(t, \theta_t, \bar{z}) = \theta_t \Pi_{C_t}(\bar{z} + \theta_t) - \tfrac{1}{2}|\bar{z} - \Pi_{C_t}(\bar{z} + \theta_t)|^2.$$

We evidently have to show that $\bar{Y}_t = \alpha Y_t$ for $t \in [0, T]$ or, equivalently, $\alpha f(t, \theta_t, \frac{z}{\alpha}) = \bar{f}(t, \theta_t, z)$. Note that for a convex set $C$, the projection $\Pi_C(a)$ is unique. If $C$ is a convex cone and $\beta > 0$, then $\beta \Pi_C(a) = \Pi_C(\beta a)$. The equality for the functions $f$ and $\bar{f}$ therefore follows. El Karoui and Rouge [7] have obtained the same BSDE and value function before Sekine.

In the following lemma we return to a technical point in the proof of Theorem 7. We show that it is possible to define a predictable process which satisfies (8). Instead of referring to a classical section theorem, see [5], we prefer to give a direct and constructive proof.

LEMMA 11 (Measurable selection). Let $(a_t)_{t \in [0,T]}$, $(\sigma_t)_{t \in [0,T]}$ be $\mathbb{R}^{1 \times m}$-valued (resp. $\mathbb{R}^{d \times m}$-valued) predictable stochastic processes, $\tilde{C} \subset \mathbb{R}^d$ a closed set and $C_t = \tilde{C}\sigma_t, t \in [0,T]$.

(a) The process

$$d = (\mathrm{dist}(a_t, \tilde{C}\sigma_t))_{t \in [0,T]}$$

is predictable.

(b) *There exists a predictable process $a^*$ with*

$$a_t^* \in \Pi_{C_t}(a_t) \qquad \text{for all } t \in [0, T].$$

PROOF. In order to prove (a), observe that $d$ is the composition of continuous mappings with predictable processes. For $k \in \mathbb{N}$, let $H^k$ denote the space of compact subsets of $\mathbb{R}^k$ equipped with the Hausdorff metric and $\mathcal{B}(H^k)$ the Borel sigma algebra with respect to this metric. The mapping $\text{dist}: \mathbb{R}^m \times H^m \to \mathbb{R}$ is jointly continuous, hence, $(\mathcal{B}(\mathbb{R}^m) \otimes \mathcal{B}(H^m) - \mathcal{B}(\mathbb{R}))$-measurable. Now consider $j: \mathbb{R}^{d \times m} \times H^d \to H^m$ that maps a compact subset $\tilde{C}$ in $\mathbb{R}^d$ by applying a $(d \times m)$-matrix $\tilde{\sigma}$ to a compact subset $\tilde{K}$ of $\mathbb{R}^m$. More formally, $j$ maps $\tilde{C}$ to the following set:

$$\tilde{K} = \{b \in \mathbb{R}^m | \exists \tilde{c} \in \tilde{C} : b = \tilde{c}\tilde{\sigma}\}.$$

The mapping $j$ is also jointly continuous and, therefore, $(\mathcal{B}(\mathbb{R}^{m \times d}) \otimes \mathcal{B}(H^d) - \mathcal{B}(H^m))$-measurable. Hence, (a) follows for compact $\tilde{C}$.

If more generally $\tilde{C}$ is closed but not bounded, take $\tilde{C}_n = \tilde{C} \cap B_n$, where $B_n$ is the closed ball with radius $n$ centered at the origin. According to what has already been shown, for $n \in \mathbb{N}$, $\text{dist}(a_t, \tilde{C}_n \sigma_t)$ defines a predictable process and $\text{dist}(a_t, \tilde{C}_n \sigma_t)$ converges to $\text{dist}(a_t, \tilde{C}\sigma_t)$, for $n \to \infty$. This proves the first claim.

In order to prove the second claim, we first concentrate on the case of compact $\tilde{C}$. We have to show that, for $z \in \mathbb{R}^m$ and a compact set $\tilde{K} \subset \mathbb{R}^m$, there exists a $(\mathcal{B}(\mathbb{R}^m) \otimes \mathcal{B}(H^m) - \mathcal{B}(\mathbb{R}^m))$-measurable mapping $\xi(z, \tilde{K})$ with $\xi(z, \tilde{K}) \in \Pi_{\tilde{K}}(z)$. This is achieved by the definition of a sequence of mappings $\xi_n(z, \tilde{K})$ with a subsequence of randomly chosen index that converges to an element of $\Pi_{\tilde{K}}(z)$. The choice of the converging subsequence will depend in a measurable way on $z$ and $\tilde{K}$.

For $n \in \mathbb{N}$, let $G_n = (x_i^n)_{i \in \mathbb{N}}$ be a dyadic grid with $\min_{x \in G_n} \text{dist}(\bar{z}, x) \leq \frac{1}{n}$ for all $\bar{z} \in \mathbb{R}^m$. Let the elements of the grid $G_n$ be numbered by $G_n = \{g_i^n : i \in \mathbb{N}\}$. Let $\tilde{K}_n$ be the elements of the grid with distance at most $\frac{1}{n}$ from $G_n$. Since we can describe the sets $\tilde{K}_n$ as the intersections of the discrete set $G_n$ with the closed set of all points in $\mathbb{R}^m$ having distance at most $\frac{1}{n}$ from $\tilde{K}$, and this closed set depends continuously on $\tilde{K}$, $\tilde{K}_n$ is measurable in $\tilde{K}$. For any $z \in \mathbb{R}^m$, let $\Pi_n(z, \tilde{K})$ be the set of all points in $\tilde{K}_n$ with minimal distance from $z$. Since $\tilde{K}_n$ is measurable in $\tilde{K}$, $\Pi_n(z, \tilde{K})$ is obviously measurable in $(z, \tilde{K})$. To define $\xi_n(z, \tilde{K})$, we have to choose one point in $\Pi_n(z, \tilde{K})$. Let it be the one with minimal index in the enumeration of $G_n$. This choice preserves the measurability in $(z, \tilde{K})$. Hence, we obtain that $\xi_n(z, \tilde{K})$ is $(\mathcal{B}(\mathbb{R}^m) \otimes \mathcal{B}(H^m) - \mathcal{B}(\mathbb{R}^m))$-measurable. Furthermore, $\liminf_{n \to \infty} |\xi_n(z, \tilde{K})| < \infty$ for all $(z, \tilde{K})$. This is one assumption in Lemma 1.55 in [8] that we aim to apply. This lemma is stated for equivalence classes of random variables, where two random variables are equivalent if they are equal almost everywhere with respect to a probability measure. Considering carefully the proof, we see that we can apply this lemma, also without reference to any measure, to obtain a result for every $(z, \tilde{K}) \in \mathbb{R}^m \times H^m$.

Lemma 1.55 in [8] yields a strictly increasing sequence $(\tau_n)_{n \in \mathbb{N}}$ of integer valued, $(\mathcal{B}(\mathbb{R}^m) \otimes \mathcal{B}(H^m) - \mathcal{B}(\mathbb{R}))$-measurable functions and a mapping $\xi : \mathbb{R}^m \times H^m \to \mathbb{R}^m$ measurable with respect to the corresponding product $\sigma$-algebra, satisfying

$$\lim_{n \to \infty} \xi_{\tau_n(z, \tilde{K})}(z, \tilde{K}) = \xi(z, \tilde{K}) \qquad \forall z \in \mathbb{R}^m, \tilde{K} \in H^m.$$

But $\xi$ is a selection. Indeed, for every $n \in \mathbb{N}$,

$$|\operatorname{dist}(z, \xi_{\tau_n}(z, \tilde{K})) - \operatorname{dist}(z, \tilde{K})| \leq \frac{1}{\tau_n} \leq \frac{1}{n}.$$

Since $\xi_{\tau_n}$ converges to $\xi$, we obtain $\operatorname{dist}(\xi, \tilde{K}) = 0$, hence, $\xi \in \tilde{K}$ and $\operatorname{dist}(z, \xi) = \operatorname{dist}(z, \tilde{K})$. Thus, by construction, $\xi(z, \tilde{K}) \in \Pi_{\tilde{K}}(z)$ for all $(z, \tilde{K}) \in \mathbb{R}^m \times H^m$.

We may then choose

$$a^* = \xi(a, C\sigma)$$

to satisfy the requirements of the second part of the assertion in the compact case.

Finally, if $\tilde{C}$ is only closed, we may proceed similarly as in the proof for (a). Let $a_t^n = \xi(a, (\tilde{C} \cap B_n)\sigma_t)$, $t \in [0, T]$. This time we apply Lemma 1.55 in [8] to the sequence of predictable processes $(a^n)_{n \in \mathbb{N}}$ and the measure $P \otimes \lambda$ on $\Omega \times [0, T]$. We obtain a strictly increasing sequence of random indices $\tilde{\tau}_n(\omega, t)$ measurable with respect to the predictable $\sigma$-algebra and a predictable process $a$ such that

$$\lim_{n \to \infty} a_t^{\tilde{\tau}_n(\omega, t)}(\omega) = a_t(\omega) \qquad \text{for } P \otimes \lambda \text{ a.e. } (\omega, t).$$

For the process $a$, we have $\operatorname{dist}(a_t, \tilde{C}\sigma_t) = 0$ $P \otimes \lambda$ a.e. $\square$

LEMMA 12. *Let $(Y, Z) \in \mathcal{H}^\infty(\mathbb{R}) \times \mathcal{H}^2(\mathbb{R}^m)$ be a solution of the BSDE (7), and let $p^*$ be given by Lemma 11 for $a = Z + \frac{1}{\alpha}\theta$. Then the processes*

$$\int_0^\cdot Z_s \, dW_s, \qquad \int_0^\cdot p_s^* \, dW_s$$

*are P-BMO martingales.*

PROOF. Let $k$ denote the upper bound of the uniformly bounded process $Y$. Applying Itô's formula to $(Y - k)^2$, we obtain, for stopping times $\tau \leq T$,

$$E\left[\int_\tau^T Z_s^2 \, ds \Big| \mathcal{F}_\tau\right] = E[(F - k)^2 | \mathcal{F}_\tau] - |Y_\tau - k|^2$$
$$- 2E\left[\int_\tau^T (Y_s - k) f(s, Z_s) \, ds \Big| \mathcal{F}_\tau\right].$$

The definition of $f$ yields, for all $(t,z) \in [0,T] \times \mathbb{R}^m$,

$$f(t,z) \leq z\theta_t + \frac{1}{2\alpha}|\theta_t|^2.$$

Therefore, there exist positive constants $c_1$, $c_2$ and $\tilde{c}_1$ such that

$$E\left[\int_\tau^T |Z_s|^2 \, ds \Big| \mathcal{F}_\tau\right] \leq c_1 + c_2 E\left[\int_\tau^T |Z_s + 1| \, ds \Big| \mathcal{F}_\tau\right]$$

$$\leq \tilde{c}_1 + \tfrac{1}{2} E\left[\int_\tau^T |Z_s|^2 \, ds \Big| \mathcal{F}_\tau\right].$$

Hence, $\int_0^\cdot Z_s \, dW_s$ is a BMO martingale.

We next deal with the stochastic integral process of $p^*$. The triangle inequality implies

$$|p^*| \leq \left|Z + \frac{1}{\alpha}\theta\right| + \left|p^* - \left(Z + \frac{1}{\alpha}\theta\right)\right|.$$

The definition of $p^*$ together with (4), yields for some constants $k_1$, $k_2$,

$$|p_t^*| \leq 2|Z_t| + \frac{2}{\alpha}|\theta_t| + k_1 \leq 2|Z_t| + k_2, \qquad t \in [0,T],$$

and, thus, for every stopping time $\tau \leq T$,

$$E\left[\int_\tau^T |p_t^*|^2 \, dt \Big| \mathcal{F}_\tau\right] \leq E\left[\int_\tau^T 8|Z_t|^2 \, dt + 2Tk_2^2 \Big| \mathcal{F}_\tau\right].$$

This implies the $P$-BMO property of $\int_0^\cdot p_s^* \, dW_s$. $\square$

**3. Power utility.** In this section we calculate the value function and characterize the optimal strategy for the utility maximization problem with respect to

$$U_\gamma(x) = \frac{1}{\gamma} x^\gamma, \qquad x \geq 0, \quad \gamma \in (0,1).$$

This time, our investor maximizes the expected utility of his wealth at time $T$ without an additional liability. The trading strategies are constrained to take values in a closed set $\bar{C}_2 \subseteq \mathbb{R}^d$. In this section we shall use a somewhat different notion of trading strategy: $\tilde{\rho} = (\tilde{\rho}^i)_{i=1,\ldots,d}$ denotes the part of the wealth invested in stock $i$. The number of shares of stock $i$ is given by $\frac{\tilde{\rho}_t^i X_t}{S_t^i}$. A $d$-dimensional $\mathbb{F}$-predictable process $\tilde{\rho} = (\tilde{\rho}_t)_{0 \leq t \leq T}$ is called trading strategy (part of wealth) if the following wealth process is well defined:

$$(11) \quad X_t^{(\tilde{\rho})} = x + \int_0^t \sum_{i=1}^d \frac{X_s^{(\tilde{\rho})} \tilde{\rho}_{i,s}}{S_{i,s}} \, dS_{i,s} = x + \int_0^t X_s^{(\tilde{\rho})} \tilde{\rho}_s \sigma_s (dW_s + \theta_s \, ds),$$

and the initial capital $x$ is positive. The wealth process $X^{(\tilde{\rho})}$ can be written as

$$X_t^{(\tilde{\rho})} = x\mathcal{E}\left(\int \tilde{\rho}_s \sigma_s (dW_s + \theta_s\, ds)\right)_t, \qquad t \in [0,T].$$

As before, it is more convenient to introduce

$$\rho_t = \tilde{\rho}_t \sigma_t, \qquad t \in [0,T].$$

Accordingly, $\rho$ is constrained to take its values in

$$C_t(\omega) = \tilde{C}\sigma_t(\omega), \qquad t \in [0,T], \omega \in \Omega.$$

The sets $C_t$ satisfy (4). In order to formulate the optimization problem, we first define the set of admissible trading strategies.

DEFINITION 13. The set of admissible trading strategies $\tilde{\mathcal{A}}$ consists of all $d$-dimensional predictable processes $\rho = (\rho_t)_{0 \le t \le T}$ that satisfy $\rho_t \in C_t(\omega) P \otimes \lambda$-a.s. and $\int_0^T |\rho_s|^2\, ds < \infty$ $P$-a.s.

Define the probability measure $Q \sim P$ by

$$\frac{dQ}{dP} = \mathcal{E}\left(-\int \theta_s\, dW_s\right)_T.$$

The set of admissible trading strategies is free of arbitrage because, for every $\rho \in \tilde{\mathcal{A}}$, the wealth process $X^{(\tilde{\rho})}$ is a local $Q$-martingale bounded from below, hence, a $Q$-supermartingale. Since $Q$ is equivalent to $P$, the set of trading strategies $\tilde{\mathcal{A}}$ is free of arbitrage.

The investor faces the maximization problem

(12) $$\bar{V}(x) = \sup_{\tilde{\rho} \in \tilde{\mathcal{A}}} E[U(X_T^{(\tilde{\rho})})].$$

In order to find the value function and an optimal strategy, we apply the same method as for the exponential utility function. We therefore have to construct a stochastic process $\tilde{R}^{(\rho)}$ with terminal value

$$\tilde{R}_T^{(\rho)} = U\left(x + \int_0^T X_s \rho_s \frac{dS_s}{S_s}\right),$$

and an initial value $\tilde{R}_0^{(\rho)} = \tilde{R}_0^x$ that does not depend on $\rho$, $\tilde{R}^{(\rho)}$ is a supermartingale for all $\rho \in \tilde{\mathcal{A}}$ and a martingale for a $\rho^* \in \tilde{\mathcal{A}}$. Then $\rho^*$ is the optimal strategy and the value function given by $\bar{V}(x) = \tilde{R}_0^x$. Applying the utility function to the wealth process yields

$$(X_t^{\rho,x})^\gamma = x^\gamma \exp\left(\int_0^t \gamma \rho_s\, dW_s + \int_0^t \gamma \rho_s \theta_s\, ds - \tfrac{1}{2}\int_0^t \gamma|\rho_s|^2\, ds\right), \qquad t \in [0,T].$$

This equation suggests the following choice:

$$\text{(13)} \quad \tilde{R}_t^{(\rho)} = x^\gamma \exp\left(\int_0^t \gamma \rho_s \, dW_s + \int_0^t \gamma \rho_s \theta_s \, ds - \tfrac{1}{2} \int_0^t \gamma |\rho_s|^2 \, ds + Y_t\right),$$

where $(Y, Z)$ is a solution of the BSDE

$$Y_t = 0 - \int_t^T Z_s \, dW_s - \int_t^T f(s, Z_s) \, ds, \qquad t \in [0, T].$$

In order to get the supermartingale property of $\tilde{R}^{(\rho)}$, we have to construct $f(t, z)$ such that, for $t \in [0, T]$,

$$\text{(14)} \quad \gamma \rho_t \theta_t - \tfrac{1}{2} \gamma |\rho_t|^2 + f(t, Z_t) \leq -\tfrac{1}{2} |\gamma \rho_t + Z_t|^2 \qquad \text{for all } \rho \in \tilde{\mathcal{A}}.$$

$\tilde{R}^{(\rho^*)}$ will even be a martingale if equality holds for $\rho^* \in \tilde{\mathcal{A}}$. This is equivalent to

$$f(t, Z_t) \leq \tfrac{1}{2} \gamma (1 - \gamma) \left|\rho_t - \frac{1}{1-\gamma}(Z_t + \theta_t)\right|^2 - \frac{1}{2} \frac{\gamma |Z_t + \theta_t|^2}{1 - \gamma} - \tfrac{1}{2} |Z_t|^2.$$

Hence, the appropriate choice for $f$ is

$$f(t, z) = \frac{\gamma(1-\gamma)}{2} \operatorname{dist}^2\left(\frac{1}{1-\gamma}(z + \theta_t), C_t\right) - \frac{\gamma|z + \theta_t|^2}{2(1-\gamma)} - \tfrac{1}{2}|z|^2,$$

and a candidate for the optimal strategy must satisfy

$$\rho_t^* \in \Pi_{C_t(\omega)}\left(\frac{1}{1-\gamma}(Z_t + \theta_t)\right), \qquad t \in [0, T].$$

In the following theorem both value function and optimal strategy are described.

THEOREM 14. *The value function of the optimization problem is given by*

$$V(x) = x^\gamma \exp(Y_0) \qquad \text{for } x > 0,$$

*where $Y_0$ is defined by the unique solution $(Y, Z) \in \mathcal{H}^\infty(\mathbb{R}) \times \mathcal{H}^2(\mathbb{R}^m)$ of the BSDE*

$$\text{(15)} \qquad Y_t = 0 - \int_t^T Z_s \, dW_s - \int_t^T f(s, Z_s) \, ds, \qquad t \in [0, T],$$

*with*

$$f(t, z) = \frac{\gamma(1-\gamma)}{2} \operatorname{dist}^2\left(\frac{1}{1-\gamma}(z + \theta_t), C_t\right) - \frac{\gamma|z + \theta_t|^2}{2(1-\gamma)} - \tfrac{1}{2}|z|^2.$$

*There exists an optimal trading strategy $\rho^* \in \tilde{\mathcal{A}}$ with the property*

$$\text{(16)} \qquad \rho_t^* \in \Pi_{C_t(\omega)}\left(\frac{1}{1-\gamma}(Z_t + \theta_t)\right).$$

PROOF. According to Lemma 11, $(f(t,z))_{t\in[0,T]}$ is a predictable stochastic process which also depends on $\sigma$. Due to (4) and the boundedness of $\theta$, Condition (H1) for Theorem 2.3 in [11] is fulfilled. We obtain the existence of a solution $(Y,Z) \in \mathcal{H}^\infty(\mathbb{R}) \times \mathcal{H}^2(\mathbb{R}^m)$ for the BSDE (15). Uniqueness follows from the comparison arguments in the uniqueness part of the proof of Theorem 7.

Let $\rho^*$ denote the predictable process constructed with Lemma 11 for $a = \frac{1}{1-\gamma}(Z+\theta)$. Lemma 17 below shows that $\rho^* \in \tilde{\mathcal{A}}$. By Theorem 2.3 in [10], the process $\tilde{R}^{(\rho^*)}$ is a martingale with terminal value

$$\tilde{R}_T^{(\rho^*)} = x^\gamma \exp\left( \int_0^T \gamma \rho_s^* \, dW_s + \int_0^T \gamma \rho_s^* \theta_s \, ds - \tfrac{1}{2} \int_0^T \gamma |\rho_s^*|^2 \, ds \right).$$

This is the power utility from terminal wealth of the trading strategy $\rho^*$. Therefore, the expected utility of $\rho^*$ is equal to $\tilde{R}_0^{(\rho^*,x)} = x^\gamma \exp(Y_0)$.

To show that this provides the value function, let $\rho \in \tilde{\mathcal{A}}$. (14) yields

$$\tilde{R}_t^{(\rho)} = x^\gamma \exp(Y_0) \mathcal{E}\left( \int (\gamma \rho_s + Z_s) \, dW_s \right)_t \exp\left( \int_0^t v_s \, ds \right), \qquad t \in [0,T],$$

for a process $v$ with $v_s \leq 0$ $\lambda \otimes P$-a.s.

The stochastic exponential is a local martingale. There exists a sequence of stopping times $(\tau_n)_{n\in\mathbb{N}}$, $\lim_{n\to\infty} \tau_n = T$ such that

$$E[\tilde{R}_{t\wedge\tau_n}^{(\rho)} | \mathcal{F}_s] \leq \tilde{R}_{s\wedge\tau_n}^{(\rho)}, \qquad s \leq t,$$

for every $n \in \mathbb{N}$. Furthermore, $\tilde{R}^{(\rho)}$ is bounded from below by 0. Passing to the limit and applying Fatou's lemma yields that $\tilde{R}^{(\rho)}$ is a supermartingale. The terminal value $\tilde{R}_T^{(\rho,x)}$ is the utility of the terminal wealth of the trading strategy $\rho$. Consequently,

$$E[U(X_T^{(\rho,x)})] \leq \tilde{R}_0^{(x)} = x^\gamma \exp(Y_0) \qquad \text{for all } \rho \in \mathcal{A}. \qquad \square$$

Again, we can show that an investor starting to act at some stopping time in the trading interval $[0,T]$ will perceive the strategy $\rho^*$ just constructed as optimal. Let $\tau \leq T$ denote a stopping time and $X_\tau$ an $\mathcal{F}_\tau$-measurable random variable which describes the capital at time $\tau$, that is, $X_\tau = X_\tau^\rho$ for a $\rho \in \tilde{\mathcal{A}}$ and an initial capital $x > 0$. Consider the maximization problem

$$(17) \qquad \bar{V}(\tau, X_\tau) = \operatorname*{ess\,sup}_{\rho \in \mathcal{A}_\tau} E\left[ U\left( X_\tau + \int_\tau^T X_s \rho_s (dW_s + \theta_s \, ds) \right) \Big| \mathcal{F}_\tau \right].$$

PROPOSITION 15 (Dynamic principle). *The value function $x^\gamma \exp(y)$ satisfies the dynamic programming principle, that is,*

$$\bar{V}(\tau, X_\tau) = (X_\tau)^\gamma \exp(Y_\tau)$$

for all stopping times $\tau \leq T$, where $Y_\tau$ is given by the unique solution $(Y, Z)$ of the BSDE (15). An optimal strategy which attains the essential supremum in (17) is given by $\rho^*$ constructed in Theorem 14.

PROOF. See Proposition 9.

REMARK 16. Suppose that the constraint set $C$ is a convex cone. Then the optimal strategy $\rho^*$ constructed in Theorem 14 is the same as in [15].

Sekine uses the utility function $x \mapsto \frac{1}{\gamma} x^\gamma$ and obtains the value function

$$\tilde{V}(x) = \frac{1}{\gamma} x^\gamma \exp((1-\gamma)\tilde{Y}_0),$$

where $\tilde{Y}_0$ is defined by the unique solution $(\tilde{Y}, \tilde{Z}) \in \mathcal{H}^\infty(\mathbb{R}) \times \mathcal{H}^2(\mathbb{R}^m)$ of the BSDE

$$\tilde{Y}_t = 0 - \int_t^T \tilde{Z}_s \, dW_s - \int_t^T g(s, \tilde{Z}_s) \, ds, \qquad t \in [0, T].$$

Here

$$g(t, \tilde{z}) = \frac{|\theta_t|^2}{2} - \frac{1}{2}\left|\theta_t - \Pi_{C_t}\left(\tilde{z} + \frac{\theta_t}{1-\gamma}\right)\right|^2 - \frac{1-\gamma}{2}\left|\tilde{z} - \Pi_{C_t}\left(\tilde{z} + \frac{\theta_t}{1-\gamma}\right)\right|^2.$$

As for the exponential utility function, we have to show $(1-\gamma)\tilde{Y} = Y$ or, equivalently, $(1-\gamma)g(t, \frac{z}{1-\gamma}) = f(t, z)$. In fact, we have

$$(1-\gamma)g\left(t, \frac{z}{1-\gamma}\right) = (1-\gamma)\left[\frac{|\theta_t|^2}{2} - \frac{1}{2}\left|\theta_t - \Pi_{C_t}\left(\frac{z+\theta_t}{1-\gamma}\right)\right|^2\right]$$

$$- \frac{(1-\gamma)^2}{2}\left|\frac{z}{1-\gamma} - \Pi_{C_t}\left(\frac{z+\theta_t}{1-\gamma}\right)\right|^2$$

$$= \theta_t \Pi_{C_t}(z+\theta_t) - \frac{1}{2(1-\gamma)}|\Pi_{C_t}(z+\theta_t)|^2$$

$$- \frac{1}{2}|z|^2 + z\Pi_{C_t}(z+\theta_t) - \frac{1}{2}|\Pi_{C_t}(z+\theta_t)|^2$$

$$= (z+\theta_t)\Pi_{C_t}(z+\theta_t) - \frac{2-\gamma}{2(1-\gamma)}|\Pi_{C_t}(z+\theta_t)|^2 - \frac{1}{2}|z|^2$$

$$= -\frac{\gamma}{2(1-\gamma)}|\Pi_{C_t}(z+\theta_t)|^2 - \frac{1}{2}|z|^2.$$

To obtain the last equality, we use

$$(z+\theta_t)\Pi_{C_t}(z+\theta_t) = |\Pi_{C_t}(z+\theta_t)|^2$$

[see (18) below].

For the function $f$, we obtain

$$\begin{aligned}f(t,z) &= \frac{\gamma(1-\gamma)}{2}\left|\frac{1}{1-\gamma}(z+\theta_t) - \Pi_{C_t}\left(\frac{1}{1-\gamma}(z+\theta_t)\right)\right|^2 \\ &\quad - \frac{\gamma}{2}\frac{(z+\theta_t)^2}{(1-\gamma)} - \frac{1}{2}|z|^2 \\ &= -\frac{\gamma}{1-\gamma}(z+\theta_t)\Pi_{C_t}(z+\theta_t) + \frac{\gamma}{2(1-\gamma)}|\Pi_{C_t}(z+\theta_t)|^2 - \frac{1}{2}|z|^2 \\ &= -\frac{\gamma}{2(1-\gamma)}|\Pi_{C_t}(z+\theta_t)|^2 - \frac{1}{2}|z|^2.\end{aligned}$$

For $t \in [0,T], z \in \mathbb{R}^m$, we therefore have

$$(1-\gamma)g\left(t, \frac{z}{1-\gamma}\right) = f(t,z).$$

It remains to prove that, for a convex cone $C$ and $a \in \mathbb{R}^m$, the following equality holds:

(18) $$\Pi_C(a)(a - \Pi_C(a)) = 0.$$

If $\Pi_C(a) = 0$, then the identity is satisfied. If not, consider the half line $\lambda \Pi_C(a)$, $\lambda \geq 0$. This half line is part of the cone $C$, so $\Pi_C(a)$ is also the projection of $a$ on the half line. □

LEMMA 17. *Let $(Y,Z) \in \mathcal{H}^\infty(\mathbb{R}) \times \mathcal{H}^2(\mathbb{R}^m)$ be a solution of the BSDE (15), and let $\rho^*$ be given by (16). Then the processes*

$$\int_0^\cdot Z_s \, dW_s, \qquad \int_0^\cdot \rho_s^* \, dW_s$$

*are P-BMO martingales.*

PROOF. We can use the same line of reasoning as in the proof of Lemma 12. The argument given there has to be slightly modified, however. We may take a lower bound $k$ for $Y$, and apply Itô's formula to $|Y-k|^2$, to conclude in the same manner as before. □

**4. Log utility.** To complete the spectrum of important utility functions, in this section we shall consider logarithmic utility. As in the preceding section, the agent has no liability at time $T$. Trading strategies and wealth process have the same meaning as in Section 3 [see (11)]. The trading strategies $\tilde{\rho}$ are constrained to take values in a closed set $\tilde{C}_2 \subset \mathbb{R}^d$. For $\rho_t = \tilde{\rho}_t \sigma_t$, the constraints are described by $C_t = \tilde{C}_2 \sigma_t, t \in [0,T]$. In order to compare the logarithmic utility of the terminal wealth of two trading strategies, we have to impose a mild integrability condition on $\rho$. Recall that $\rho^i > 1$ means that the investor has to borrow money in order to buy stock $i$ and if $\rho^i < 0$, then the investor has a negative number of stock $i$. An integrability condition on $\rho$ is not restrictive.

DEFINITION 18. The set of admissible trading strategies $\mathcal{A}_l$ consists of all $\mathbb{R}^d$-valued predictable processes $\rho$ satisfying $E[\int_0^T |\rho_s|^2 \, ds] < \infty$ and $\rho_t \in C_t$ $P \otimes \lambda$-a.s.

For the logarithmic utility function,

$$U(x) = \log(x), \qquad x > 0,$$

we obtain a particularly simple BSDE that leads to the value function and the optimal strategy. The optimization problem is given by

$$
\begin{aligned}
V(x) &= \sup_{\rho \in \mathcal{A}_l} E[\log(X_T^{(\rho)})] \\
&= \log(x) + \sup_{\rho \in \mathcal{A}_l} E\left[\int_0^T \rho_s \, dW_s + \int_0^T (\rho_s \theta_s - \tfrac{1}{2}|\rho_s|^2) \, ds\right],
\end{aligned}
\tag{19}
$$

where the initial capital $x$ is positive again. As in Section 2, we want to determine a process $R^{(\rho)}$ with $R_T^{(\rho)} = \log(X_T^{(\rho)})$, and an initial value that does not depend on $\rho$. Furthermore, $R^{(\rho)}$ is a supermartingale for all $\rho \in \mathcal{A}_l$, and there exists a $\rho^* \in \mathcal{A}_l$ such that $R^{(\rho^*)}$ is a martingale. The strategy $\rho^*$ is the optimal strategy and $R_0^{\rho^*}$ is the value function of the optimization problem (19).

We can choose, for $t \in [0, T]$,

$$R_t^{(\rho)} = \log x + Y_0 + \int_0^t (\rho_s + Z_s) \, dW_s + \int_0^t (-\tfrac{1}{2}|\rho_s - \theta_s|^2 + \tfrac{1}{2}\theta_s^2 + f(s)) \, ds,$$

where

$$f(t) = \tfrac{1}{2} \operatorname{dist}^2(\theta_t, C_t) - \tfrac{1}{2}|\theta_t|^2, \qquad t \in [0, T],$$

and $(Y_t, Z_t)$ is the unique solution of the following BSDE:

$$Y_t = 0 - \int_t^T Z_s \, dW_s - \int_t^T f(s) \, ds, \qquad t \in [0, T].$$

Due to Definition 18, the boundedness of $\theta$ and (4), the stochastic integral in $R^{(\rho)}$ is a martingale for all $\rho \in \mathcal{A}_l$. Hence, $R^{(\rho)}$ is a supermartingale for all $\rho \in \mathcal{A}_l$. An optimal trading strategy $\rho^*$, which satisfies $\rho_t^* \in \Pi_{C_t}(\theta_t)$, can be constructed by means of Lemma 11. The initial value $Y_0$ satisfies

$$Y_0 = -E\left[\int_0^T f(s) \, ds\right].$$

Hence,

$$V(x) = R_0^{\rho^*}(x) = \log(x) + E\bigg[-\int_0^T f(s)\,ds\bigg].$$

In particular, $\rho^*$ only depends on $\theta$, $\sigma$ and the set $\tilde{C}_2$ describing the constraints on the trading strategies.

**Acknowledgments.** The authors wish to thank the two anonymous referees and an Associate Editor for their helpful comments and critics. Our paper has improved greatly during the revision process because of that.

Y. Hu
IRMAR
Campus de Beaulieu
Université de Rennes 1
F-35042 Rennes Cedex
France
e-mail: ying.hu@univ-rennes1.fr

P. Imkeller
M. Müller
Institut für Mathematik
Humboldt-Universität zu Berlin
Unter den Linden 6
D-10099 Berlin
Germany
e-mail: imkeller@mathematik.hu-berlin.de
e-mail: muellerm@mathematik.hu-berlin.de